\documentclass[11pt,reqno]{elsarticle}

\usepackage{placeins}

\makeatletter
\def\ps@pprintTitle{%
 \let\@oddhead\@empty
 \let\@evenhead\@empty
 \def\@oddfoot{}%
 \let\@evenfoot\@oddfoot}
\makeatother

\usepackage{amsmath,amsthm,amscd,wrapfig,amsfonts,amssymb,mathtools}
\usepackage{enumitem,comment,mathabx,mathrsfs}
\usepackage{booktabs,multirow,lscape,cool}
\usepackage{url,parskip,caption}
\usepackage[top=1.0in, bottom=1.0in, left=1.0in, right=1.0in]{geometry}

\newtheorem{remark}{Remark}

\newtheorem{conjecture}{Conjecture}

\providecommand{\e}[1]{\ensuremath{\times 10^{#1}}}

\makeatletter
\g@addto@macro\normalsize{%
  \setlength\abovedisplayskip{.4em}
  \setlength\belowdisplayskip{.4em}
  \setlength\abovedisplayshortskip{.4em}
  \setlength\belowdisplayshortskip{.4em}
}

\providecommand{\PS}[1]{\ensuremath{\mathit{Ps}\hskip.05em{}_{#1}}}
\providecommand{\QS}[1]{\ensuremath{\mathit{Qs}\hskip.05em{}_{#1}}}

\providecommand{\PsiP}[1]{\ensuremath{\mathit{\Psi S}\hskip.05em{}_{#1}}}

\begin{document}

\begin{frontmatter}


\title
{
An $\mathcal{O}\left(1\right)$ algorithm  for the numerical evaluation of the Sturm-Liouville eigenvalues 
of the spheroidal wave functions of order zero
}

\begin{abstract}
In addition to being the eigenfunctions of the restricted Fourier operator,
the angular spheroidal wave functions of the first kind of order zero and nonnegative 
integer characteristic exponents are the solutions of a singular self-adjoint 
Sturm-Liouville problem.    The running time of the standard algorithm for the 
numerical evaluation of their Sturm-Liouville eigenvalues grows with both bandlimit 
and characteristic exponent.  Here, we describe a new approach 
whose running time is bounded independent of these parameters.
Although the Sturm-Liouville eigenvalues
are of little interest themselves, our algorithm is a component
of a fast scheme for the numerical evaluation of the
prolate spheroidal wave functions  developed by one of the  authors.   
We illustrate the performance of our method with numerical experiments.

\end{abstract}

\author[1]{Rafeh Rehan}
\author[1]{James Bremer\corref{cor1}}
\ead{bremer@math.toronto.edu}
\cortext[cor1]{Corresponding author}

\address[1]{Department of Mathematics,  University of Toronto}


\end{frontmatter}
\begin{section}{Introduction}

The angular prolate spheroidal wave functions of the first kind of order zero and nonnegative integer characteristic exponents
\begin{equation}
\PS{0}(z;\gamma),\ \PS{1}(z;\gamma),\ \PS{2}(z;\gamma),\ \ldots
\label{introduction:ps}
\end{equation}
%
are the eigenfunctions of the restricted Fourier operator
\begin{equation}
\mathscr{F}_\gamma\left[f\right](z) = \int_{-1}^{1} \exp(i\gamma t z) f(t)\ dt.
\label{introduction:fourier}
\end{equation}
 As such they provide an efficient  mechanism for representing 
functions in the image of $\mathscr{F}_\gamma$,
which is the space of functions with bandlimit $\gamma$.   
Indeed,  the magnitudes of the first $2/\pi \gamma$ eigenvalues of the restricted Fourier operator
are close to  $\sqrt{2\pi/\gamma}$,
the magnitudes of the next $\mathcal{O}\left(\log(\gamma)\right)$ eigenvalues
decay extremely rapidly, and the remaining eigenvalues are all close
to zero  \cite{Landau-Widom}.  It follows that only the first
$2/\pi \gamma + \mathcal{O}\left(\log(\gamma)\right)$ functions in (\ref{introduction:ps})
are needed to represent elements of  the image of $\mathscr{F}_\gamma$
with high relative accuracy.

The behaviour of the spectrum of the restricted Fourier operator
makes the numerical calculation of $\PS{n}(z;\gamma)$ and the corresponding 
eigenvalue $\lambda_n(\gamma)$ through the direct discretization
of  (\ref{introduction:fourier}) extremely difficult.
Fortunately, the functions (\ref{introduction:ps})
are also the solutions of the singular self-adjoint Sturm-Liouville problem
\begin{equation}
\left\{
\begin{aligned}
(1-z^2)y''(z) -2z y'(z) + (\chi-\gamma^2z^2) y(z) &= 0, \ \ \ -1 < z < 1, \\
\lim_{z \to  \pm 1} y'(z) \sqrt{1-z^2}&=0
\end{aligned}
\right.
\label{introduction:slp}
\end{equation}
(see, for instance, Section~3.8 of \cite{Meixner}).  We refer to the differential equation
in (\ref{introduction:slp}) as the reduced spheroidal wave equation because it is obtained 
from the more familiar spheroidal wave equation by deleting one of its parameters (order).

The Osipov-Xiao-Rokhlin method 
\cite{Xiao-Rokhlin-Yarvin,Osipov-Rokhlin-Xiao}, which is the standard approach to the numerical
calculation of $\PS{n}(z;\gamma)$ and the corresponding
Sturm-Liouville eigenvalue $\chi_n(\gamma)$,
operates by representing a solution  of  (\ref{introduction:slp})
as a finite Legendre expansion.    While the dependence of its running time 
on the parameters $\gamma$ and $n$ is not fully understood,
the numerical experiments of \cite{Feichtinger} suggest that it grows
as $\mathcal{O}\left(n + \sqrt{n \gamma}\right)$, at least for large values
of $n$ and $\gamma$.

In \cite{prolates1}, a numerical scheme for calculating $\PS{n}(z;\gamma)$ which runs in time
 independent of $n$ and which grows sublogarithmically with $\gamma$ is described.   However, it requires
knowledge of the value of $\chi_n(\gamma)$.
Here, we describe a mechanism for evaluating $\chi_n(\gamma)$ with near machine precision accuracy
in time independent of $\gamma$ and $n$.  
It proceeds by constructing a piecewise polynomial expansion of  a nonstandard  
analytic continuation $\chi_\xi(\gamma)$ of $\chi_n(\gamma)$.
The parameter $\xi$ is related to the  value of a certain phase function for the reduced
spheroidal wave equation at the point $0$, and we also construct expansions
which allow for the rapid evaluation of the values of the first few derivatives of this 
phase function at $0$.    The ability to rapidly evaluate these quantities allows us 
to accelerate the algorithm of \cite{prolates1}, reducing its running time
by a factor of 10 or so.

Many second order differential equations admit phase functions
which are easier to represent using standard mechanisms (such as polynomial expansions)
than the solutions of the equations themselves.  This is often demonstrated by
proving that the equation admits a modulus function which satisfies various monotonicity
properties.  It is well known that Legendre's differential equation, which is 
a special case of the reduced spheroidal wave equation, possesses such a modulus
function (see, for instance, \cite{durand78}).  Here, we conjecture that
the reduced spheroidal wave equation admits a modulus function with properties similar
to this modulus function for  Legendre's differential equation.  
We also present the results of numerical experiments showing that, in any case,
$\chi_\xi(\gamma)$ can be represented extremely efficiently via polynomial expansions.
Indeed, the expansion of $\chi_\xi(\gamma)$ we constructed for this article
 consumes less than $0.76$ MB of memory and 
allows for evaluation of $\chi_n(\gamma)$
for all $2^{6} \leq \gamma \leq 2^{20}$ and $0 \leq n \leq 1.1 \gamma$.
Each evaluation takes less than $5 \times 10^{-6}$ seconds on the standard desktop computer
used to conduct the experiments of this paper.  The range $0 \leq n \leq 1.1 \gamma$ was chosen because
\begin{equation}
\lambda_{\lfloor 1.1 \gamma\rfloor} (\gamma) < \epsilon_0
\end{equation}
for all $\gamma \geq 2^6 = 64$, where  $\epsilon_0 =2^{-52} \approx 2.220446049250313 \times 10^{-16}$ is machine zero for 
IEEE double precision arithmetic.    The Osipov-Xiao-Rokhlin algorithm is more efficient than
the approach suggested here for values of $\gamma$ smaller than $64$, and it is to be preferred in that
regime.
However, expansions which hold for a larger range of the parameters, including smaller values of 
$\gamma$, could easily be constructed.

The properties of $\chi_\xi(\gamma)$ are in stark contrast to those of the 
standard analytic continuation $\chi_\nu(\gamma)$ of $\chi_n(\gamma)$
obtained via characteristic exponents (see, for instance, \cite{Meixner}).
The latter  is entire in $\gamma$, but only meromorphic in $\nu$, with branch points at each half-integer
value of $\nu$.   This greatly complicates any attempt to construct expansions
of $\chi_\nu(\gamma)$ using standard machinery, like polynomial or trigonometric
expansions.

The remainder of this article is structured as follows.  
Section~\ref{section:phase} briefly discusses phase functions for second order 
linear ordinary differential equations.     In Section~\ref{section:spheroidal}, we 
define certain standard solutions of the
reduced spheroidal wave equation and a define a particular phase function
for the reduced spheroidal wave equation which plays a central role in our algorithm.
Section~\ref{section:spheroidal} includes a discussion of characteristic exponents and the standard analytic continuations of 
$\chi_n(\gamma)$ and the spheroidal wave functions.
In Section~\ref{section:monotonicity}, we give several conjectures regarding the properties
of a particular phase function for the reduced spheroidal wave equation and discuss some consequences
of these conjectures.  
In Section~\ref{section:index}, we introduce an alternative to characteristic exponents for indexing
the reduced spheroidal wave functions.
Section~\ref{section:algorithm} details our numerical algorithm.  
In Section~\ref{section:experiments}, we present the results of numerical experiments demonstrating the 
properties of our algorithm.   

\end{section}

\begin{section}{Phase functions for second order differential equations}
\label{section:phase}

Suppose that $\Omega$ is a simply-connected open set in the complex plane, and
that $q:\Omega\to\mathbb{C}$ is an analytic function.
Then we say that $\psi:\Omega \to\mathbb{C}$ is a phase function for the second
order linear ordinary differential equation
\begin{equation}
y''(z) + q(z) y(z) = 0,\ \ \ z \in \Omega,
\label{phase:ode}
\end{equation}
provided $\psi'$ does not vanish on $\Omega$ and 
\begin{equation}
u(z) = \frac{\sin(\psi(z))}{\sqrt{\psi'(z)}} \ \ \mbox{and}\ \ 
v(z) = \frac{\cos(\psi(z))}{\sqrt{\psi'(z)}}
\label{phase:uv}
\end{equation}
form a basis in the space of solutions of (\ref{phase:ode}).  The particular realization of the square root
used in (\ref{phase:uv}) is immaterial.  It can be verified through a straightforward calculation
that $\psi'$ satisfies the second order nonlinear ordinary differential equation
\begin{equation}
q(z) - (\psi'(z))^2 + \frac{3}{4} \left(\frac{\psi''(z)}{\psi'(z)}\right)^2
-\frac{1}{2}\frac{\psi'''(z)}{\psi'(z)}=0,
\label{phase:kummer}
\end{equation}
which we call Kummer's equation after E.E.~Kummer who studied it in \cite{Kummer}.  Conversely,
if $\psi'$ does not vanish in $\Omega$ and satisfies (\ref{phase:kummer}) then the function $u$ and $v$
defined via (\ref{phase:uv}) are solutions of (\ref{phase:ode}).
In light of (\ref{phase:uv}), we refer to
\begin{equation}
m(z) = \frac{1}{\psi'(z)} = (u(z))^2 + (v(z))^2
\label{phase:m}
\end{equation}
as the modulus function associated with the phase function $\psi'$.

If $u, v$ is a pair of solutions of (\ref{phase:ode}) whose (necessarily constant) Wronskian $w$ is nonzero on $\Omega$
and such that the modulus function (\ref{phase:m}) does not vanish on $\Omega$, 
then it can be easily verified that the function
\begin{equation}
\psi'(z) = \frac{w}{(u(z))^2 + (v(z))^2}
\label{phase:psip}
\end{equation}
satisfies Kummer's equation.  It follows that  any antiderivative $\psi$ of $\psi'$ is a phase function for
(\ref{phase:ode}).  Adding the requirement that (\ref{phase:uv}) holds determines $\psi$ up to an additive
constant multiple of $2\pi$.



\end{section}

\begin{section}{The prolate spheroidal wave functions of order zero}
\label{section:spheroidal}

In this section, we discuss characteristic exponents, review the definitions of some of the standard
solutions of the spheroidal wave function and define a certain phase function for 
the reduced spheroidal wave equation.

\begin{subsection}{The spheroidal wave equation}

The spheroidal wave equation
\begin{equation}
(1-z^2) y''(z) -2z y'(z) + \left( \chi - \gamma^2 z^2 - \frac{\mu^2}{1-z^2} \right) y(z) = 0 
\label{spheroidal:swe}
\end{equation}
arises when the method of separation of variables is used to solve the constant coefficient Helmholtz equation
(see, for instance, Chapter~5 of \cite{Morse-FeshbachI}).
When $\gamma^2 >0$, its solutions are known as prolate spheroidal wave functions,
and when $\gamma^2<0$ they are known as the oblate spheroidal wave functions.
The spheroidal wave functions are typically indexed via the explicit parameters $\gamma$ and $\mu$, 
which we refer to as the bandlimit and order, respectively, 
and by an implicit parameter $\nu$  known as the characteristic exponent.
The explicit parameter $\chi$ in (\ref{spheroidal:swe}) is usually regarded
as a function of $\gamma$, $\mu$ and $\nu$.

In this article, we restrict our attention to the prolate spheroidal wave functions
of order zero (i.e., we impose the restrictions $\gamma^2>0$ and  $\mu=0$).  These are the spheroidal wave
functions which are the most widely used in applications.  Obviously, they
are solutions of the differential equation
\begin{equation}
(1-z^2) y''(z) -2z y'(z) + \left( \chi - \gamma^2 z^2 \right) y(z) = 0,
\label{spheroidal:rswe}
\end{equation}
which we call the reduced spheroidal wave equation.   It has regular singular points at $z = \pm 1$ 
and an irregular singular point at infinity.

\end{subsection}

\begin{subsection}{Characteristic exponents}

For any complex value of the parameter $\chi$,  (\ref{spheroidal:rswe}) admits a solution of the form  
\begin{equation}
z^\nu \sum_{n=-\infty}^\infty a_n z^{2n}
\label{spheroidal:char1}
\end{equation}
with the Laurent expansion convergent in the annulus $1 < |z| < \infty$
and, in the event that $\nu$ is not a half-integer,
there is a second solution of the form
\begin{equation}
z^{-\nu-1} \sum_{n=-\infty}^\infty b_n z^{2n},
\label{spheroidal:char2}
\end{equation}
also with the  Laurent expansion convergent in the annulus $1 < |z| < \infty$
(see, for instance, \cite{Ince} or \cite{Hille}).
The parameter $\chi$  appearing in (\ref{spheroidal:rswe}) obviously only determines the value of $\nu$ 
up to an integral multiple of $2$.  It is less obvious
that for each value of $\nu$ which is not a half-integer, there 
are a countable collection of values of $\chi$ such that (\ref{spheroidal:rswe}) admits
a pair of solutions, one of which is of the form (\ref{spheroidal:char1})
and the other is of the form (\ref{spheroidal:char2}).  A proof of this appears in \cite{Meixner}.

It is standard (again, see \cite{Meixner}) to associate  a unique value of  $\chi$, which we denote by $\chi_\nu(\gamma)$,
with each $\nu$ which is not a half-integer by requiring that
\begin{equation}
\lim_{\gamma\to0^+} \chi_\nu(\gamma) =  \nu(\nu+1).
\label{spheroidal:chilimit}
\end{equation}
This condition is motivated by that fact that  (\ref{spheroidal:rswe}) reduces to Legendre's differential equation when $\gamma=0$. 
In that case, there is a solution of the form
\begin{equation}
z^\nu \sum_{n=0}^\infty a_n z^{2n},
\end{equation}
where $a_0 \neq 0$, and it can be easily seen that $\chi$ relates to $\nu$ via the formula $\chi = \nu(\nu+1)$.

The function $\chi_\nu(\gamma)$ which results from imposing the condition (\ref{spheroidal:chilimit})
is analytic in $\gamma$, but only meromorphic in $\nu$, with branch points at the half-integers.
Figure~\ref{figure:chiplot}, which contains a plot of $\chi_\nu(\gamma)$ as a function of $\nu$ when $\gamma=2$,
shows the jump discontinuities that occur at  half-integer values of $\nu$.
Moreover,  because $\nu(\nu+1) = (-\nu-1)(-\nu)$, we have that $\chi_\nu(\gamma) = \chi_{-\nu-1}(\gamma)$.
However, since the value of $\chi$ determines the possible values of $\nu$ in (\ref{spheroidal:char1}) and (\ref{spheroidal:char2})
up to an additive constant which is an integral multiple of $2$,
these are the only two values of $\nu$ which can correspond to a particular choice of $\chi$.

While it is not possible to extend this scheme in order to assign a unique value of $\chi$ to each half-integer value of $\nu$,
one can associate two distinct values of $\chi$ to each half-integer $\nu$ by taking 
limits from the left and right.  Further information on the case of half-integer characteristic exponents
can be found in \cite{Imam}.  In what follows, we will assume implicitly that $\nu$ is not a half-integer and 
this will cause no difficulties for us.

\begin{figure}[!h]
\begin{center}
\includegraphics[width=.6\textwidth]{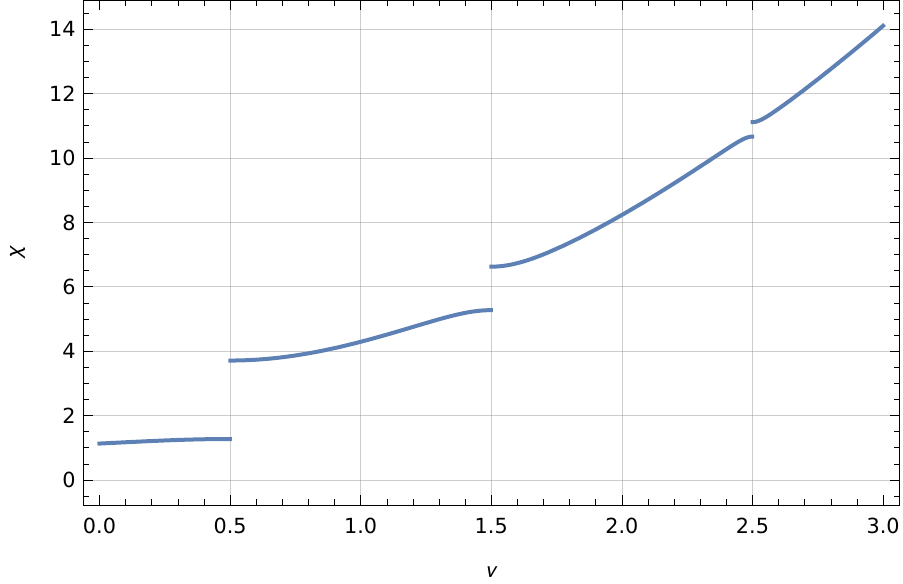}
\end{center}
\caption{A plot of the standard analytic continuation $\chi_\nu(\gamma)$ of $\chi_n(\gamma)$ as a function of $\nu$ when $\gamma=2$.  There is a jump
discontinuity at each half-integer value of $\nu$. }
\label{figure:chiplot}
\end{figure}

\end{subsection}

\begin{subsection}{The angular prolate spheroidal wave functions of the first and second kinds of order zero}

Zero is a double root of the indicial equation for the reduced spheroidal wave equation corresponding to the regular singular point at $z=1$.
Accordingly,  it admits a one-dimensional subspace of solutions which are regular at $1$
and a one-dimensional subspace of solutions which have logarithmic singularities at $1$ (see, for instance, \cite{Hille}).

We use  $\PS{\nu}(z;\gamma)$ to denote the unique solution of (\ref{spheroidal:rswe})  which is regular at $z=1$ and such that
either the value of $\PS{\nu}(0;\gamma)$ agrees with that of the Legendre function $P_\nu(0)$ or, in the event
that  $P_\nu(0)=0$, the derivative of  $\PS{\nu}(x;\gamma)$ with respect to $x$ at $0$ agrees with the derivative of the Legendre function
$P_\nu(x)$ at $0$.  We refer to $\PS{\nu}(z;\gamma)$ as the angular spheroidal wave function of the first kind of bandlimit $\gamma$, order zero
and characteristic exponent $\nu$.  It is well known that $\PS{\nu}(z;\gamma)$  admits an expansion of the form
\begin{equation}
\PS{\nu}(z;\gamma) = \sum_{k=-\infty}^\infty a_k (\nu;\gamma) P_{\nu+2k}(z)
\label{spheroidal:ps}
\end{equation}
and, like the Legendre functions of the first kind,
$\PS{\nu}(z;\gamma)$ is entire as a function of $z$ when $\nu$ is an integer
and is single-valued on the cut-plane $\mathbb{C}\setminus\left(-\infty,-1\right]$
when $\nu$ is not an integer. 

Assuming $\nu$ is not an integer, it can be easily verified that 
\begin{equation}
\QS{\nu}(z;\gamma) = \sum_{k=-\infty}^\infty a_k (\nu;\gamma) Q_{\nu+2k}(z),
\label{spheroidal:qs}
\end{equation}
where $Q_\nu$ is the Legendre function of the second kind of degree $\nu$ and the coefficients $\{a_k(\nu;\gamma)\}$ are the same as
in (\ref{spheroidal:ps}), is also solution of the reduced spheroidal wave equation.
The representation (\ref{spheroidal:qs}) is problematic when $\nu$ is an integer because $Q_\nu(z)$, when viewed as a function of $\nu$, 
has simple poles at the negative integers.  However, in this case, it is possible to find a representation of the form
\begin{equation}
\QS{\nu}(z;\gamma) = \sum_{k=-\infty}^\infty b_k (\nu;\gamma) \frac{Q_{\nu+2k}(z)}{\Gamma(\nu+2k+1)}
\label{spheroidal:qs2}
\end{equation}
since  $\Gamma(\nu+1)$ has simple zeros at each negative integer (see, for instance, \cite{Meixner}).
We refer to $\QS{\nu}(z;\gamma)$ as the angular prolate spheroidal wave function  of the second kind of bandlimit $\gamma$,
order zero and characteristic exponent $\nu$.  Just like the Legendre functions of the second kind,
the function $\QS{\nu}(z;\gamma)$ is defined for $z$ on the cut-plane $\mathbb{C}\setminus\left(-\infty,1\right]$
and has a logarithmic singularity at the point $z=1$.

\begin{remark}
Some of the formulas in this article become simpler when $\PS{\nu}(z,\gamma)$ is normalized 
through the requirement that $\PS{\nu}(1,\gamma)=1$.  This is in keeping with the
standard convention for the normalization of the Legendre functions of the first kind.  However, since many
of the angular spheroidal wave functions decay exponentially on some portion
of the interval $(0,1)$, such a normalization scheme would result
in some of the $\PS{\nu}(z,\gamma)$  taking on extremely large values 
in the interval $(0,1)$, thus complicating their numerical evaluation.
\end{remark}

\end{subsection}

\begin{subsection}{The angular prolate spheroidal wave functions of the first kind of order zero and integer characteristic exponents }
The boundary conditions
\begin{equation}
\lim_{z \to \pm 1} y'(z) (1-z^2) = 0
\label{spheroidal:bc}
\end{equation}
together with Equation~(\ref{spheroidal:rswe}) comprise a singular self-adjoint Sturm-Liouville problem
(see, for example, \cite{Zettl} for a discussion of such problems.).
The angular prolate spheroidal wave functions of the first kind of order zero and nonnegative integer characteristic exponents
\begin{equation}
\PS{0}(z;\gamma),\ \PS{1}(z;\gamma),\ \PS{2}(z;\gamma),\ \ldots
\label{spheroidal:eigenfunctions}
\end{equation}
are a collection of solutions of this Sturm-Liouville problem which form an orthogonal basis in $L^2(-1,1)$.
Much of the interest in the spheroidal wave functions
stems from the fact that (\ref{spheroidal:eigenfunctions}) are also
eigenfunctions of the restricted Fourier operator
\begin{equation}
\mathscr{F}_\gamma\left[f\right](z) = \int_{-1}^1 \exp\left(i\gamma z t \right) f(t)\ dt.
\end{equation}
This observation was widely publicized in the article \cite{PSWFI} published in the 1960s,
but it was known much earlier (see, for instance, Section~3.8 of \cite{Meixner} and the references cited there).

\end{subsection}

\begin{subsection}{The radial prolate spheroidal wave functions of order zero}
\label{subsection:s3}
Another solution of (\ref{spheroidal:rswe}), which is known as the radial spheroidal wave function
of the third kind of order zero, is given by the formula
\begin{equation}
S^{(3)}_\nu(z;\gamma) =  \frac{1}{\PS{\nu}(1,\gamma)}\int_1^\infty \exp(i\gamma zt) \PS{\nu}(t; \gamma)\ dt.
\label{spheroidal:s3}
\end{equation}
The integral is absolutely convergent for $\mbox{Im}(z) > 0$ and
$S^{(3)}_\nu(z;\gamma)$ is typically taken to be its analytic continuation
to the cut plane $\mathbb{C}\setminus\left(-\infty,1\right]$.
However, it is more convenient for us to regard the domain of 
$S^{(3)}_\nu(z;\gamma)$ as the analytic continuation of (\ref{spheroidal:s3}) to 
an open simply-connected set $\Gamma$ containing $\{ z : \mbox{Im}(z) \geq 0\ \mbox{and}  \ z \neq \pm 1\}$.
The asymptotic behaviour of $S^{(3)}_\nu(z;\gamma)$ can be easily deduced from (\ref{spheroidal:s3}):
\begin{equation}
S^{(3)}(z;\gamma) = \frac{\exp\left(i \gamma z\right)}{\gamma z} + \mathcal{O}\left(\frac{1}{z^2}\right)  \ \ \mbox{as} \ \ z \to\infty.
\label{spheroidal:s3asym}
\end{equation}

Similarly, the radial prolate spheroidal wave function of the fourth kind
of bandlimit $\gamma$, order zero and characteristic exponent $\nu$ is given by the
formula
\begin{equation}
S^{(4)}_\nu(z;\gamma) = 
\frac{1}{\PS{\nu}(1;\gamma)}
\int_1^\infty \exp(-i\gamma zt) \PS{\nu}(t; \gamma)\ dt.
\label{spheroidal:s4}
\end{equation}
The integral is absolutely convergent for $\mbox{Im}(z) < 0$, and, like 
$S^{(3)}_\nu(z;\gamma)$, the domain of this function is usually taken to be 
the cut plane $\mathbb{C}\setminus \left(-\infty,1\right]$.  However,
we regard $S^{(4)}_\nu(z;\gamma)$ as defined on the same open simply-connected
set $\Gamma$ which serves as the domain of $S^{(3)}_\nu(z;\gamma)$.
It follows easily from (\ref{spheroidal:s4}) that
\begin{equation}
S^{(4)}(z;\gamma) = \frac{\exp\left(-i \gamma z\right)}{\gamma z} + \mathcal{O}\left(\frac{1}{z^2}\right)  \ \ \mbox{as} \ \ z \to\infty.
\label{spheroidal:s4asym}
\end{equation}

We define the radial prolate spheroidal wave functions of the first and second kinds
of bandlimit $\gamma$, order zero and characteristic exponent $\nu$
on $\Gamma$ via the formulas
\begin{equation}
S^{(1)}_\nu(z;\gamma) = \frac{S^{(3)}_\nu(z;\gamma) + S^{(4)}_\nu(z;\gamma)}{2}
\label{spheroidal:s1}
\end{equation}
and
\begin{equation}
S^{(2)}_\nu(z;\gamma) = \frac{S^{(3)}_\nu(z;\gamma) - S^{(4)}_\nu(z;\gamma)}{2i},
\label{spheroidal:s2}
\end{equation}
respectively.  Then
\begin{equation}
S^{(3)}_\nu(z;\gamma) =  S^{(1)}_\nu(z;\gamma) +  iS^{(2)}_\nu(z;\gamma).
\end{equation}
and
\begin{equation}
S^{(4)}_\nu(z;\gamma) =  S^{(1)}_\nu(z;\gamma) -  iS^{(2)}_\nu(z;\gamma).
\end{equation}
Moreover, from (\ref{spheroidal:s3asym}) and (\ref{spheroidal:s4asym})
it easily follows that 
%
\begin{equation}
S^{(1)}(z;\gamma) = \frac{\sin\left(\gamma z\right)}{\gamma z} + \mathcal{O}\left(\frac{1}{z^2}\right)  \ \ \mbox{as} \ \ z \to\infty
\label{spheroidal:s1asym}
\end{equation}
and
\begin{equation}
S^{(2)}(z;\gamma) = \frac{\cos\left(\gamma z\right)}{\gamma z} + \mathcal{O}\left(\frac{1}{z^2}\right)  \ \ \mbox{as} \ \ z \to\infty.
\label{spheroidal:s2asym}
\end{equation}
%



\label{section:spheroidal:radial}
\end{subsection}

\begin{subsection}{The normal form of the reduced spheroidal wave equation}

A straightforward calculation shows that if $y$ satisfies the reduced spheroidal wave equation (\ref{spheroidal:rswe}), then
the function $u(z) = y(z) \sqrt{1-z^2}$ satisfies
\begin{equation}
u''(z) + \left(\frac{1}{\left(1-z^2\right)^2} + \frac{\chi-\gamma^2z^2}{1-z^2}\right) u(z) = 0, \ \ z \in \Gamma.
\label{spheroidal:nrswe}
\end{equation}
The particular realization of $\sqrt{1-z^2}$ used here is immaterial.  We refer to (\ref{spheroidal:nrswe}) as the normal
form of the reduced spheroidal wave equation.

\end{subsection}

\begin{subsection}{The phase and modulus functions associated with $S^{(3)}_\nu(z;\gamma)$}

We define the functions $M_\nu(z;\gamma)$ and $\PsiP{\nu}(z;\gamma)$  on $\Gamma$ via the formulas
\begin{equation}
M_\nu(z;\gamma) =  \left(S^{(1)}_\nu(z;\gamma)\right)^2 +  \left(S^{(2)}_\nu(z;\gamma)\right)^2
\end{equation}
and
\begin{equation}
\PsiP{\nu}(z;\gamma) = \int_1^z \frac{\gamma}{M_\nu(u;\gamma) (1-u^2)}\ du.
\label{spheroidal:psip}
\end{equation}
Since the Wronskian of any pair of solutions of the differential equation (\ref{spheroidal:nrswe}) is constant,
it can be easily seen from (\ref{spheroidal:s1asym}) and (\ref{spheroidal:s2asym}) that the Wronskian of 
the pair $S^{(1)}_\nu(z;\gamma)\sqrt{1-z^2},\ S^{(2)}_\nu(z;\gamma)\sqrt{1-z^2}$ is $\gamma$.
It follows that $\PsiP{\nu}(z;\gamma)$
is a phase function for the normal form of the reduced spheroidal wave equation, and that $M_\nu(z;\gamma)(1-z^2)$ 
is the corresponding modulus function.  We omit the factor of $(1-z^2)$ in the definition of $M_\nu(z;\gamma)$ 
to make stating the conjectures of Section~\ref{section:monotonicity} more convenient.  Moreover, by a slight abuse
of terminology, we will refer to $M_\nu(z;\gamma)$ as the modulus function associated with $\PsiP{\nu}(z;\gamma)$.

Since $\PsiP{\nu}(z;\gamma)$ is a phase function for (\ref{spheroidal:nrswe}), 
there exist  $C_\nu(z;\gamma)$ and $D_\nu(z;\gamma)$ such that
\begin{equation}
\begin{aligned}
\PS{\nu}(z;\gamma)\sqrt{1-z^2} &= 
C_\nu(z;\gamma)  \frac{\sin\left(\PsiP{\nu}(z;\gamma)\right)}{\sqrt{\frac{d\PsiP{\nu}}{dz}(z,\gamma)}}
+  D_\nu(z;\gamma)
\frac{\cos\left(\PsiP{\nu}(z;\gamma)\right)}{\sqrt{\frac{d\PsiP{\nu}}{dz}(z,\gamma)}}.
\end{aligned}
\label{spheroidal:ps1}
\end{equation}
From (\ref{spheroidal:psip}), we see that (\ref{spheroidal:ps1}) is equivalent to 
\begin{equation}
\PS{\nu}(z;\gamma)  = 
C_\nu(z;\gamma)  \frac{\sqrt{M_\nu(z;\gamma)}}{\sqrt{\gamma}} \sin\left(\PsiP{\nu}(z;\gamma)\right)
+ D_\nu(z;\gamma)  \frac{\sqrt{M_\nu(z;\gamma)}}{\sqrt{\gamma}} \cos\left(\PsiP{\nu}(z;\gamma)\right).
\label{spheroidal:ps2}
\end{equation}
Because $\QS{\nu}(z;\gamma)$ has a logarithmic singularity at $1$, we must have
\begin{equation}
\lim_{z\to 1} \left|\sqrt{M_\nu(z;\gamma)}\right| = \infty.
\label{spheroidal:limit1}
\end{equation}
But we also have
\begin{equation}
\lim_{z\to 1}\PsiP{\nu}(z;\gamma) = 0,
\end{equation}
and it follows from this and (\ref{spheroidal:limit1}) that
\begin{equation}
\lim_{z \to 1} 
\left|
\sqrt{M_\nu(z;\gamma)}\cos\left(\PsiP{\nu}(z;\gamma)\right)
\right| = \infty.
\end{equation}
Since $\PS{\nu}(z;\gamma)$ is nonsingular at $1$, we must have $D_\nu(\gamma) =0 $ in (\ref{spheroidal:ps1}) and 
(\ref{spheroidal:ps2}) so that
\begin{equation}
\PS{\nu}(z;\gamma)\sqrt{1-z^2} = C_\nu(z;\gamma) \frac{\sin\left(\PsiP{\nu}(z;\gamma)\right)}{\sqrt{\frac{d\PsiP{\nu}}{dz}(z,\gamma)}}.
\label{spheroidal:psc}
\end{equation}

\end{subsection}

\begin{subsection}{The reduced spheroidal wave functions as functions of the parameter $\chi$}

It follows easily both from mechanism used to define  $\PS{\nu}(z;\gamma)$ and from 
(\ref{spheroidal:ps}) that $\PS{-\nu-1}(z;\gamma) = \PS{\nu}(z;\gamma)$.  
In particular,  $\PS{\nu}(z;\gamma)$ is uniquely determined by the value of  the parameter $\chi$ in (\ref{spheroidal:rswe}). 
From this observation and the definitions of Section~\ref{section:spheroidal:radial},
it is clear that the radial spheroidal wave functions, and hence also $M_\nu(z;\gamma)$ and $\PsiP{\nu}(z;\gamma)$,
are  uniquely determined by $\chi$ and hence can be indexed via $\chi$ instead of by $\nu$.  We note that
it follows from (\ref{spheroidal:qs}) that  this is not the case for $\QS{\nu}(z;\gamma)$.

We will, by a slight abuse of notation,  use $S^{(3)}_\chi(z;\gamma)$ to 
denote the radial spheroidal wave functions of the third kind corresponding to 
 $\chi = \chi_\nu(\gamma)$, and likewise for $M_\chi(z;\gamma)$ and $\PsiP{\chi}(z;\gamma)$.
It follows from standard results in the theory of ordinary differential equations
that these functions are entire in $\chi$ as well as in $\gamma$.  

\label{section:spheroidal:chi}
\end{subsection}



%

\end{section}

\begin{section}{The monotonicity properties of the reduced spheroidal wave equation}
\label{section:monotonicity}

It is well known that many second order differential equations admit modulus functions
which satisfy strong monotonicity properties.  Bessel's differential equation
furnishes one such  example.      The formula
\begin{equation}
J_\lambda^2(z) + Y_\lambda^2(z) = 
 \frac{2}{\pi}\int_0^\infty \exp(-zt) P_{\lambda-\frac{1}{2}} \left(1+\frac{t^2}{2}\right)\ dt,
\label{monotonicity:nicholson}
\end{equation}
which can be found in \cite{hartman73},  expresses 
a modulus function for Bessel's equation as the Laplace transform of a positive function.
Because of the close relationship between modulus and phase functions, it follows that Bessel's equation
admits a phase functions which is, among other things, increasing and nonoscillatory on the interval $(0,\infty)$.
This is in stark contrast to the Bessel functions themselves, which behave as increasing or decreasing
exponential functions on the interval $\left(0,\sqrt{\lambda^2-1/4}\right)$ and oscillate on 
$\left(\sqrt{\lambda^2-1/4},\infty\right)$.  The existence of this phase function was used at an early
date to rapidly evaluate the Bessel functions \cite{Goldstein-Thaler} of large arguments,  and it is 
exploited by the widely used algorithm \cite{Amos} for the same purpose.

Similar results hold for many second order linear ordinary differential equations.
 Relevant formulas  for the Jacobi functions, Gegenbauer functions and Hermite functions can be found in \cite{durand78},
and the articles \cite{hartman61} and \cite{hartman73} 
give conditions under which a second order linear ordinary differential equation
admits a modulus function which is the Laplace transform of a nonnegative Borel measure.

The asymptotic estimates (\ref{spheroidal:s3asym}) and  (\ref{spheroidal:s4asym}) indicates that, at least for large $z$,
the modulus function $M_\nu(z;\gamma)$ can be well approximated by a nonoscillatory function.
This suggests that $M_\nu(z;\gamma)$, like (\ref{monotonicity:nicholson}), satisfies various monotonicity
properties.  Our suspicions are further bolstered by the fact that 
Legendre's differential equation, which is a special case of (\ref{spheroidal:rswe}),
is known to satisfy certain strong monotonicity properties (see \cite{durand78}).
In this section,  after briefly defining various notions of monotonicity, we make several conjectures
about the monotonicity properties of $M_\nu(z;\gamma)$ and $S_\nu(z;\gamma)$.
These conjectures were arrived at through numerical experiments, experiments using computer algebra systems
and our belief that the properties of the reduced spheroidal wave equation are similar to
those of Legendre's differential equation.

\begin{subsection}{Monotonicity Properties}
A smooth function $f$ defined on an open interval $I$ is completely monotone if
$(-1)^j f^{(j)}(z) \geq 0$ for all nonnegative integers $j$ and all $z \in I$.  
It is absolutely monotone provided
provided $f^{(j)}(z) \geq 0$ for all nonnegative integers $j$ and all $z \in I$.
A $k$-times differentiable function $f$ defined on an open interval $I$ is 
$k$-times monotone provided $(-1)^j f^{(j)}(z) \geq 0$ for all nonnegative integers $j \leq k$ and all $z \in I$.

It is well known  that $f$ is completely monotone
on $(0,\infty)$ if and only if $f$ is the Laplace transform of a nonnegative Borel
measure (see, for instance, \cite{Widder}).
Similarly, a function $f$ is $k$-times monotone, where $k \geq 1$,
on $(0,\infty)$ if and only if there is a nonnegative Borel measure $\alpha$ such that 
\begin{equation}
f(x) = \int_0^{\frac{1}{x}} (1-tx)^{(k-1)}\ d\alpha(t).
\label{monotonicity:ktimes}
\end{equation}
Formula (\ref{monotonicity:ktimes}) suggests an obvious generalization of the notion of $k$-times monotone to noninteger values of $k$;
that is, a function $f$ is said to be $\omega$-times monotone, where $\omega \geq 1$ is not necessarily an integer, provided
there is a nonnegative Borel measure $\alpha$ such that
\begin{equation}
f(x) = \int_0^{\frac{1}{x}} (1-tx)^{(\omega-1)}\ d\alpha(t).
\label{monotonicity:omegatimes}
\end{equation}
See \cite{williamson} for a proof of (\ref{monotonicity:ktimes}) and a discussion of the definition (\ref{monotonicity:omegatimes}).
\end{subsection}


\begin{subsection}{Conjectures regarding  $S^{(3)}_\nu(z;\gamma)$, $M_\nu(z;\gamma)$  and $\PsiP{\nu}(z;\gamma)$}

We now state several conjectures regarding $S^{(3)}_\nu(z;\gamma)$ and the associated  phase and modulus functions.
The first of these asserts that the properties the modulus function $M_\nu(z;\gamma)$ mirror those
of a modulus function for Legendre's differential equation.

\begin{conjecture}
For fixed $\gamma > 0$ and $\nu > 0$, when viewed as a function of $z$,
$M_\nu(z;\gamma)$  is absolutely monotone on $(0,1)$ and completely monotone on $(1,\infty)$.
\label{conjecture1}
\end{conjecture}

The second of our conjectures can be viewed as a stronger version of the Sturm Comparison Theorem
in that it implies that the reduced spheroidal wave functions become
more oscillatory as the parameter $\chi$ increases.
\begin{conjecture}
For fixed $\gamma > 0$ and $-1 <z <1$, 
$\PsiP{\chi}(z;\gamma)$ is strictly decreasing on the interval $(0,\infty)$.
\label{conjecture2}
\end{conjecture}
%

%
%

Finally, we have following conjecture which generalizes one made in \cite{prolates1} to the case
of noninteger values of $\nu$:
\begin{conjecture}
For a fixed $\gamma > 0$ and $\nu > 0$, when viewed as a function of $z$,
$S^{(3)}_{\nu}(iz;\gamma)$ is $(\nu+2)$-times monotone on $(0,\infty)$.
\label{conjecture4}
\end{conjecture}

\end{subsection}

\end{section}

\begin{section}{An alternative method for indexing the reduced spheroidal wave functions}
\label{section:index}
Characteristic exponents are the standard mechanism for indexing the reduced spheroidal wave functions.
This scheme has advantage that the solutions of most interest --- those which are the eigenfunctions
of the restricted Fourier operator --- correspond to nonnegative integer characteristic exponents.
However, it has the serious disadvantage that $\chi_\nu(\gamma)$ has branch points at the
half-integer values of $\nu$.  

In Section~\ref{section:spheroidal:chi}, we observed that many of the spheroidal wave functions,
as well as the phase and modulus functions defined in this article,  can be indexed via the parameter $\chi$
appearing in the reduced spheroidal wave equation.  It follows from standard results in the theory of ordinary
differential equations that these functions are entire in both $\gamma$ and $\chi$.  However, 
the values of $\chi$ corresponding to the nonnegative integer characteristic exponents
are not apparent and require substantial effort to calculate.

We now discuss a new mechanism for indexing the reduced spheroidal wave functions which
combines the advantages of both of these approaches.
It is a consequence of (\ref{spheroidal:psc})  that the zeros of $\PS{\nu}(z;\gamma)$ 
must occur at points $z$ such that $\PsiP{\nu}(z;\gamma)$ is an integral multiple of $\pi$.
It is well known that when $n$ is a nonnegative integer, $\PS{n}(z;\gamma)$ has $n$ zeros in the interval $[0,1)$,
that it has a zero at the point $z=0$ if $n$ is odd, and that its first derivative has a zero
at $z=0$ if $n$ is even (see, for instance, \cite{Meixner}).     It follows from these observations
that
\begin{equation}
\PsiP{n}(0;\gamma) = -\frac{\pi}{2} \left(n+1\right)
\label{spheroidal:psn}
\end{equation}
whenever $n$ is an nonnegative integer.  Moreover, it is a consequence
of Conjecture~\ref{conjecture2} that the map $\chi \to -\PsiP{\chi}(0;\gamma)$ can be inverted.
 This suggests that we use the new parameter
\begin{equation}
\xi = -\frac{2}{\pi} \PsiP{\chi}(0;\gamma) - 1
\label{index:xi}
\end{equation}
to index the reduced spheroidal wave functions.   From (\ref{spheroidal:psn}), we see that,
just like characteristic exponents, nonnegative integer values of $\xi$ correspond to  the eigenfunctions 
of the restricted Fourier operator.        Indeed, the parameter $\xi$ 
generalizes the notion of ``the number of zeros of the function $\PS{\nu}(z;\gamma)$ on 
the interval $[0,1]$.''    
Because the qualitative behaviour of the reduced spheroidal wave function $\PS{n}(z;\gamma)$ is related
to the  ratio of the characteristic exponent $n$ to $\gamma$, we find it slightly more convenient
to use the parameter
\begin{equation*}
\sigma = \frac{\xi}{\gamma}
\end{equation*}
to index $\chi$ and the spheroidal wave functions.
 We denote the value of $\chi$ corresponding to $\sigma$ and $\gamma$ by $\chi_\sigma(\gamma)$.

Figure~\ref{figure:chiplot2} contains  plots of $\chi_\sigma(\gamma)$ as a  
function of $\sigma$ for several values of $\gamma$.  We note that each of these graphs
have inflection points when $\sigma \approx 2/\pi$.  
There is a regime change when $\sigma$ is somewhat larger than $2/\pi$.  For smaller
values of $\sigma$,  the reduced spheroidal wave equation has turning points in the interval $(0,1)$.
But for larger values of $\sigma$,  the spheroidal wave functions
of order zero are oscillatory on the entire interval $(-1,1)$.

\begin{figure}[!h]
\begin{center}

\includegraphics[width=.45\textwidth]{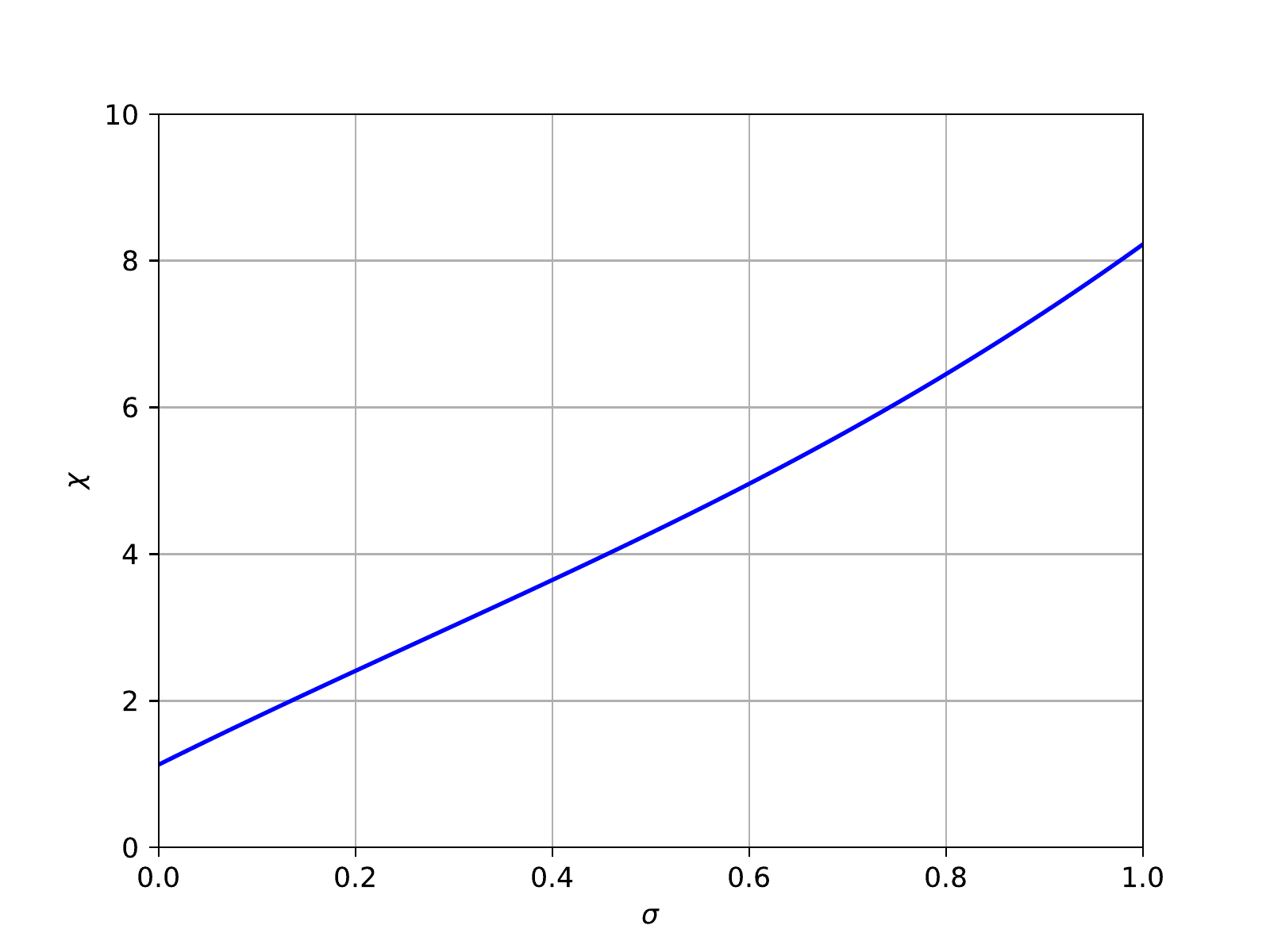}
\hfil
\includegraphics[width=.45\textwidth]{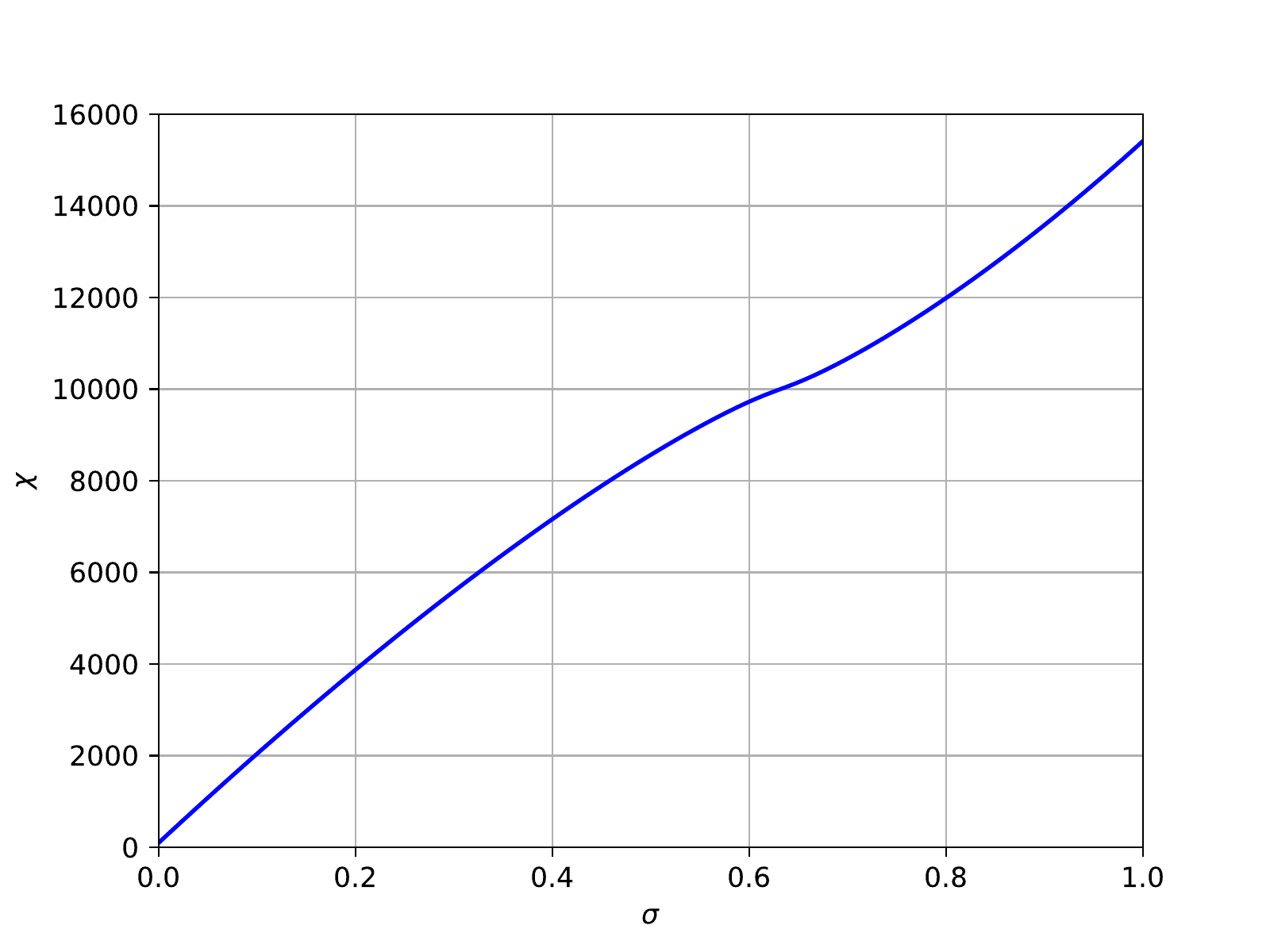}

\includegraphics[width=.45\textwidth]{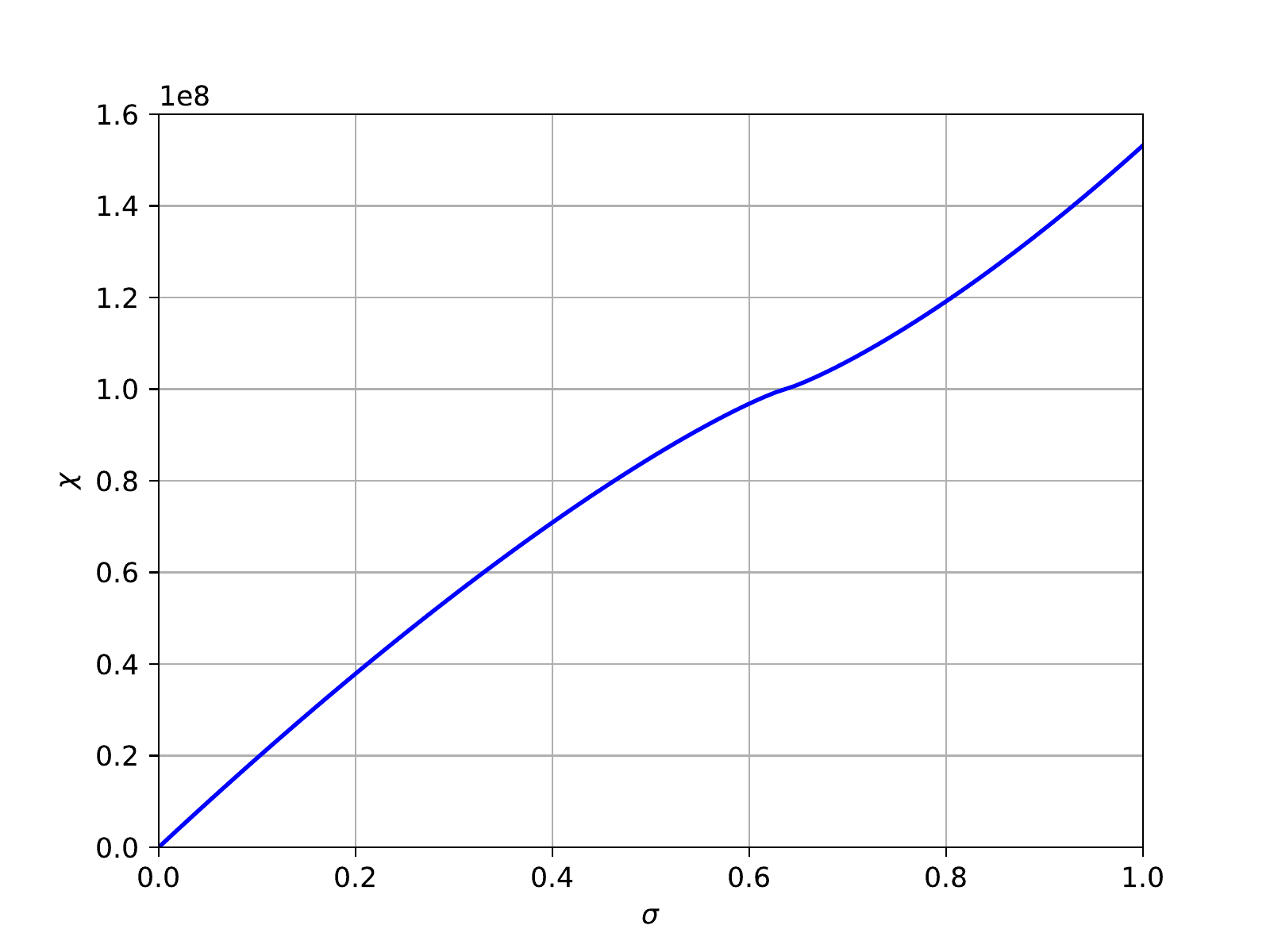}
\hfil
\includegraphics[width=.45\textwidth]{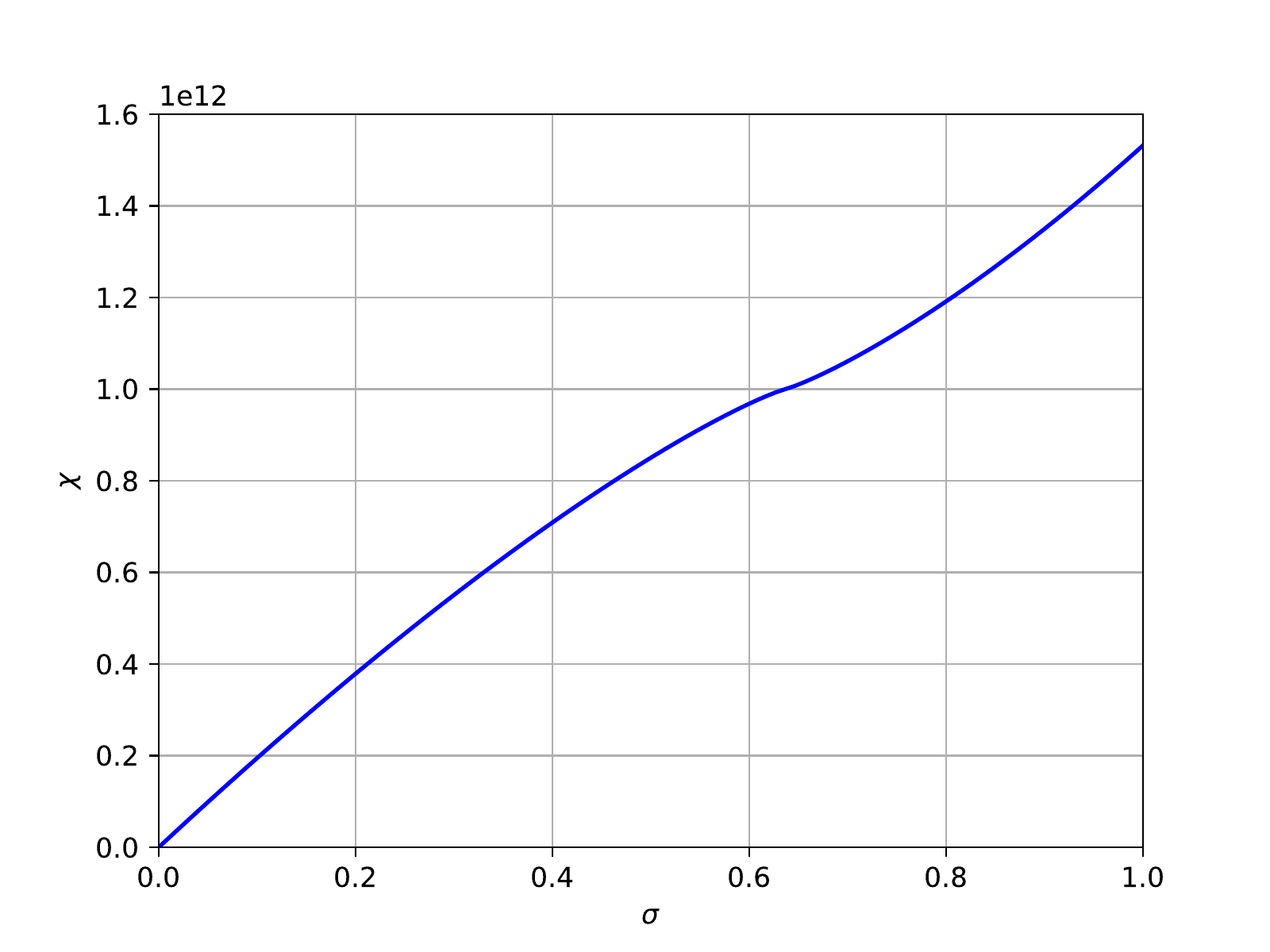}

\end{center}
\caption{Plots of the parameter $\chi$ as a function of $\sigma$ when $\gamma = 2$ (upper left), $\gamma=100$ (upper right),
$\gamma=10,000$ (lower left) and when $\gamma=10^6$ (lower right).  The relationship between $\xi$
and $\sigma$ is $\xi = \gamma \sigma$.}
\label{figure:chiplot2}
\end{figure}

\end{section}

\begin{section}{Numerical algorithm}
\label{section:algorithm}

We now describe our method for the numerical evaluation of the Sturm-Liouville eigenvalues
of the reduced spheroidal wave equation.  In broad outline, it consists of precomputing a
piecewise bivariate polynomial expansion of the function $\chi_\sigma(\gamma)$
which can then be used to rapidly evaluate $\chi_n(\gamma)$ for 
all $2^6 \leq \gamma \leq 2^{20}$
and all  nonnegative integer values of  $n$ in the interval $[0,1.1 \gamma]$.
For values of $\gamma$ smaller than $2^6$, the
Osipov-Xiao-Rokhlin algorithm is more efficient and should be preferred.

Knowledge of $\chi_n(\gamma)$ is required by the algorithm of \cite{prolates1} 
for the rapid evaluation of $\PS{n}(z;\gamma)$.  Moreover, its running time
  can be substantially accelerated when the values of 
\begin{equation}
\frac{d \PsiP{n}}{dz} (0;\gamma),\ 
\frac{d^2 \PsiP{n}}{dz^2} (0;\gamma)
 \ \mbox{and} \ 
\frac{d^3 \PsiP{n}}{dz^3} (0;\gamma)
\label{algorithm:psiders}
\end{equation}
are known.  The technique described here to construct an expansion of 
$\chi_\sigma(\gamma)$ was also used to construct expansions of the quantities (\ref{algorithm:psiders}).
However, the method was so similar to  the procedure used to construct the expansion of $\chi_\sigma(\gamma)$
that we omit the details.

We first describe the form of the expansion used to represent  $\chi_\sigma(\gamma)$.  
Then, we describe the method used to construct it. 

\begin{subsection}{The mechanism used to represent $\chi_\sigma(\gamma)$}

A $k$-term Chebyshev expansion on the interval $[a,b]$
is a sum of the form
\begin{equation}
f(x) = \sum_{j=0}^{k-1} a_j T_j\left(\frac{2}{b-a} x + \frac{a+b}{a-b} \right),
\label{algorithm:chebexp}
\end{equation}
where $T_j$ denotes the Chebyshev polynomial of degree $j$.  It is well known 
that the coefficients in the expansion  (\ref{algorithm:chebexp}) can be evaluated in a numerically stable fashion
given its values at the $k$ points
\begin{equation}
\frac{b+a}{2} - \frac{b-a}{2} \cos\left(\frac{j}{k} \pi \right), \ \ \ j=0,\ldots,k-1.
\label{algorithm:chebnodes}
\end{equation}
The set  (\ref{algorithm:chebnodes}) is known as  the $k$-point Chebyshev extrema grid on the interval $[a,b]$.
Moreover, the barycentric Lagrange formula can be used to evaluate (\ref{algorithm:chebexp})
in a numerically stable fashion  at any point in $[a,b]$ given its values at the points (\ref{algorithm:chebnodes}).
See, for instance, \cite{Trefethen} for a thorough discussion of Chebyshev interpolation.

A piecewise $k$-term Chebyshev expansion comprises a partition
\begin{equation}
a = x_0 < x_1 < x_2 < \ldots < x_n = b
\label{algorithm:part}
\end{equation}
together with a collection of $k$-term Chebyshev expansions,
one for each of the subintervals $[x_j,x_{j+1}]$.  As with Chebyshev expansions, a piecewise
Chebyshev expansion can be evaluated in a numerically stable fashion given either the coefficients
in each of these expansions, or the values of each expansion at the nodes of the 
$k$-point Chebyshev extrema grid on the  corresponding interval.  

For each $l=1,\ldots,7$ we use $I_l$ to denote the interval
\begin{equation}
I_l = \left[4^{2+l}, 4^{3+l}\right].
\label{algorithm:ik}
\end{equation}
Moreover, for each $l=1,\ldots,7$, we let
\begin{equation}
\gamma_0^{(l)},\ldots,\gamma_{k-1}^{(l)}
\end{equation}
be the nodes of the $k$-point Chebyshev extrema grid on the interval $I_l$,
where $k=30$. For each $l=1,\ldots,7$ and each  $i=0,\ldots,(k-1)$  we use a piecewise Chebyshev expansion
on the interval $(0,1.1)$  to represent  the function
\begin{equation}
f_i^{(l)}(\sigma) = \chi_\sigma\left(\gamma_i^{(l)}\right).
\label{algorithm:expos}
\end{equation}
The number of terms in each of these Chebyshev expansions is $k=30$, but the associated partitions of $(0,1.1)$ vary.
We found experimentally that the intervals $I_l$ were a suitable partition of the domain
of $\gamma$.  The partitions for the piecewise Chebyshev expansions of the functions
(\ref{algorithm:expos}) were  determined via an adaptive algorithm which is described 
in the next subsection.

Using the piecewise Chebyshev expansions of the functions (\ref{algorithm:expos}), 
$\chi_\sigma(\gamma)$ can be evaluated for all $64 = 2^6 \leq \gamma \leq 2^{20} = 1,048,576$ and 
all $0 \leq \sigma \leq 1.1 \gamma$.   More explicitly, given a pair 
of parameters $\gamma$ and $\sigma$ at which we wish to evaluate the expansion,
we first find an interval $I_l$ containing $\gamma$ ($\gamma$ might be on the boundary between
two of the subintervals, in which case either subinterval will serve).     Next we evaluate
\begin{equation}
\chi_\sigma\left(\gamma_0^{(l)}\right), \ldots, \chi_\sigma\left(\gamma_{k-1}^{(l)}\right)
\label{algorithm:chivals}
\end{equation}
using the piecewise Chebyshev expansions of the functions (\ref{algorithm:expos}).
Finally, we use the barycentric Lagrange interpolation formula for the Chebyshev 
polynomials to evaluate $\chi_\sigma(\gamma)$ using the values (\ref{algorithm:chivals}).

\end{subsection}

\begin{subsection}{Construction of the expansion}

In \cite{prolates1}, an algorithm for calculating
the phase function $\PsiP{\chi}(z;\gamma)$ for given values of the parameters $\gamma$
and $\chi$ is described.  Here, we detail how it can be used to construct a piecewise polynomial
expansion of the function
\begin{equation}
f(\sigma) = \chi_\sigma\left(\gamma\right)
\end{equation}
over the interval $0 \leq \sigma \leq 1.1$ for a fixed value of $\gamma$.
This technique is, of course, applied with $\gamma$ taking on each of the values
$$\gamma_i^{(l)}, \ \ l=1,\ldots,7, , i=0,\ldots,29.$$
in order to construct the piecewise polynomial expansions of the functions (\ref{algorithm:expos}).

As a first step, we calculate the value $\tilde{\chi}_1$ of the parameter $\chi$ corresponding
to the characteristic exponent $0$ and the value $\tilde{\chi}_2$ of the parameter $\chi$ corresponding
to the characteristic exponent  $m = \lceil 1.1 \gamma\rceil$ 
using the Osipov-Xiao-Rokhlin method.  Next, we construct a $k$-term piecewise Chebyshev
expansion --- with $k$ again taken to be  $30$ ---  which represents the function
\begin{equation}
g(\chi) = -\frac{2}{\pi} \PsiP{\chi}(0;\gamma) -1
\end{equation}
over the interval $[\tilde{\chi}_1,\tilde{\chi}_2]$.
We do this via an adaptive algorithm which operates as follows.  It maintains two lists of intervals,
one a list of processed intervals and the other a list of intervals to process.  Initially,
the interval $[\tilde{\chi}_1,\tilde{\chi}_2]$ is in the list of intervals to process and the list of
processed intervals is empty.  As long as the list of intervals to process is not empty, 
the following procedure is  repeated:

\begin{enumerate}
\item Remove an interval $[a,b]$ from the list of intervals to process.
\item Use the algorithm of \cite{prolates1} to evaluate the function $g(\chi)$ at each of the nodes
in the $k$-point Chebyshev extrema grid on $[a,b]$.
\item Compute the coefficients $a_0,\ldots,a_{k-1}$ in the Chebyshev expansion
\begin{equation*}
 \sum_{j=0}^{k-1} a_j T_j\left(\frac{2}{b-a} x + \frac{a+b}{a-b} \right)
\end{equation*}
which is equal to $g(\chi)$ at each of the nodes in the $k$-point Chebyshev extrema
grid on $[a,b]$.

\item If 
\begin{equation*}
\sum_{j=k/2}^{k-1} a_j^2 < 100\ \epsilon_0^2 \sum_{j=0}^{k-1} a_j^2,
\end{equation*}
where $\epsilon_0$ is machine zero, then
move the interval $[a,b]$ into the list of processed intervals.  Otherwise, add the intervals
$[a,(a+b)/2]$ and $[(a+b)/2,b]$ to the list of intervals to process.

\end{enumerate}
When the process terminates, the list of processed intervals determines the partition of $[\tilde{\chi}_1,\tilde{\chi}_2]$
used by the piecewise Chebyshev expansion of $g(\chi)$.

The next step consists of constructing a piecewise Chebyshev expansion
for the inverse function $f(\sigma)$ of $g(\chi)$.   
This is accomplished via an algorithm which is quite similar to that used 
to construct the expansion of $g(\chi)$.  It also maintains a list of intervals
to be processed and a list of processed intervals.  
Initially, $[0,1.1]$ is placed
in the list of intervals to process and the list of processed intervals is empty.
The algorithm then repeats the following steps until the list of intervals to process is empty:
\begin{enumerate}
\item 
Remove an interval $[c,d]$ from the list of intervals to process.

\item
Find $\alpha$ and $\beta$ such that $g(\alpha) \leq c < d \leq g(\beta)$.
The existence of $\alpha$ and $\beta$ with these properties follows from
the choice of the interval $[\tilde{\chi}_1,\tilde{\chi}_2]$ over which we represent the function
$g$.  Moreover, $\alpha$ and $\beta$ can be found by examining the values of $g(\chi)$ at the nodes
of the $k$-point Chebyshev grids on each of the subintervals associated with the piecewise Chebyshev
expansion of $g$.

\item
For each node $\sigma_j$ in the $k$-point Chebyshev extrema grid on the interval $[a,b]$
compute the value of $\chi_j$ such that $g(\chi_j)=\sigma_j$ via bisection.  Of course,
we will also have $f(\sigma_j) = \chi_j$.

More explicitly, we repeat the following steps until the quantity $\left|\alpha-\beta\right|/\left|\alpha\right|$ falls below
$10 \epsilon_0$:
\begin{itemize}
\item Use the piecewise Chebyshev expansion of $g(\chi)$ to compute $g\left(\frac{\alpha+\beta}{2}\right)$.
\item If $g\left(\frac{\alpha+\beta}{2}\right) < \sigma_j$, then let $\alpha = \frac{\alpha+\beta}{2}$.
\item Otherwise, let $\beta = \frac{\alpha+\beta}{2}$.
\end{itemize}

\item
Form the coefficients $b_0,\ldots,b_{k-1}$ in the Chebyshev expansion
\begin{equation*}
\sum_{j=0}^{k-1} b_j T_j\left(\frac{2}{d-c} x + \frac{c+d}{c-d} \right)
\end{equation*}
which, for each $j=0,\ldots,k-1$, takes on the value $\chi_j$ at the point $\sigma_j$.

\item If 
\begin{equation*}
\sum_{j=k/2}^{k-1} b_j^2 < 100\ \epsilon_0^2 \sum_{j=0}^{k-1} b_j^2,
\end{equation*}
where $\epsilon_0$ is machine zero, then
move the interval $[c,d]$ into the list of processed intervals.  Otherwise, add the intervals
$[c,(c+d)/2]$ and $[(c+d)/2,d]$ to the list of intervals to process.

\end{enumerate}

The list of processed intervals determines the partition of $[0,1.1]$ used by the piecewise Chebyshev expansion
of the inverse function $f(\sigma)$ of $g(\chi)$.

\end{subsection}

\end{section}

\begin{section}{Numerical experiments}
\label{section:experiments}

In this section, we present the results of numerical experiments which were conducted
to  to illustrate the effectiveness of the algorithm of this article.
The code for these experiments was written in Fortran and compiled with version 11.1.0 of the the GNU 
Fortran compiler.  They were performed on a desktop computer equipped with an AMD Ryzen 3900X processor.
An implementation of our algorithm and  code for conducting all of the experiments
discussed here is available on GitHub at the following address:
\begin{center}
\url{https://github.com/JamesCBremerJr/Prolates}
\end{center}

The expansion of $\chi_\sigma(\gamma)$ used in these experiments allows us to evaluate
it over the following ranges of the parameters:
\begin{equation*}
2^6 \leq \gamma \leq 2^{20} \ \ \mbox{and} \ \ 0 \leq \sigma \leq 1.1.
\end{equation*}
It occupies less than $0.76$ megabytes of memory.

In some of these experiments, we compared the performance of our algorithm
with that of the Osipov-Xiao-Rokhlin method \cite{Osipov-Rokhlin-Xiao}.  Its running time is highly
dependent on the dimension of the tridiagonal matrix formed in order to calculate
$\chi_n(\gamma)$.   Most implementations  use a highly conservative value for this dimension.  
The authors of \cite{Osipov-Rokhlin-Xiao}, for instance, take it to be  $1000 + n + \left \lfloor 1.1 \gamma \right\rfloor$ 
in their implementation.  The experiments of  \cite{Feichtinger}, however, suggest that the necessary
dimension grows as  $\mathcal{O}\left(n + \sqrt{n\gamma}\right)$.
It is difficult to find a simple formula which suffices in all cases
of interest.  Accordingly,  our implementation of the Osipov-Rokhlin-Xiao algorithm initially
takes the estimate to be 
\begin{equation}
50 + \left\lfloor\frac{2}{\pi} n\right\rfloor + \left\lfloor \sqrt{\gamma n}\right\rfloor,
\end{equation}
which we found to be sufficient for a large range of parameters, and 
then increases the dimension adaptively as needed to ensure high accuracy.
Our implementation can be found in the GitHub repository mentioned above.

To account for the vagaries of modern computational environments, all times reported here 
were obtained by repeating each calculation $100$ times and averaging the
result.


\begin{subsection}{The accuracy with which $\chi_n(\gamma)$ is calculated}

In this first set of experiments, we measured the accuracy with which $\chi_\sigma(\gamma)$
is calculated for various ranges of values of the parameters.  In each experiment,
we fixed a range $[\gamma_1,\gamma_2]$ of values of $\gamma$ and a range $[\sigma_1,\sigma_2]$ 
of values of $\sigma$.  We sampled $100$ random values of $\gamma$ in  $[\gamma_1,\gamma_2]$ and $100$ random 
integers $n$ in the interval $[\gamma \sigma_1, \gamma \sigma_2]$.
For each of the $10,000$ pairs of the sampled parameters, we calculated 
$\chi_n(\gamma)$ and compared the result with that obtained by running
the Osipov-Xiao-Rokhlin algorithm using extended precision (Fortran REAL*10) arithmetic.
Our algorithm was executed using double precision arithmetic.  We used extended precision
for the Osipov-Xiao-Rokhlin algorithm  because it looses a few digits of 
accuracy for certain values of the parameters, mainly in cases in which  $\gamma$ is large 
and $n$ is small.  Table~\ref{experiments:table1} reports the results.    Each row there
reports the maximum relative error in $\chi_n(\gamma)$ encountered for
each range of values of the parameters considered.

\begin{table}[!h]
\begin{center}
\small
\begin{tabular}{ccr@{\hspace{2em}}ccr}
\toprule
Range of $\gamma$ & Range of $\sigma$ & Max relative error  &Range of $\gamma$ & Range of $\sigma$ & Max relative error  \\\midrule
$4^3$ to $4^4$                                                                                       &$0.00 - 0.25 $ & 4.95\e{-15} &$4^7$ to $4^8$                                                                                       &$0.00 - 0.25 $ & 4.92\e{-15}   \\ & 
$0.25 - 0.50 $ & 4.82\e{-15} & & 
$0.25 - 0.50 $ & 3.91\e{-15}   \\ & 
$0.50 - 0.75 $ & 5.61\e{-15} & & 
$0.50 - 0.75 $ & 4.06\e{-15}   \\ & 
$0.75 - 1.00 $ & 5.33\e{-15} & & 
$0.75 - 1.00 $ & 4.38\e{-15}   \\\addlinespace[.25em]
$4^4$ to $4^5$                                                                                       &$0.00 - 0.25 $ & 5.05\e{-15} &$4^8$ to $4^9$                                                                                       &$0.00 - 0.25 $ & 5.12\e{-15}   \\ & 
$0.25 - 0.50 $ & 4.03\e{-15} & & 
$0.25 - 0.50 $ & 4.35\e{-15}   \\ & 
$0.50 - 0.75 $ & 4.80\e{-15} & & 
$0.50 - 0.75 $ & 4.18\e{-15}   \\ & 
$0.75 - 1.00 $ & 5.11\e{-15} & & 
$0.75 - 1.00 $ & 4.44\e{-15}   \\\addlinespace[.25em]
$4^5$ to $4^6$                                                                                       &$0.00 - 0.25 $ & 4.85\e{-15} &$4^9$ to $4^{10}$                                                                                    &$0.00 - 0.25 $ & 4.68\e{-15}   \\ & 
$0.25 - 0.50 $ & 4.13\e{-15} & & 
$0.25 - 0.50 $ & 3.82\e{-15}   \\ & 
$0.50 - 0.75 $ & 4.27\e{-15} & & 
$0.50 - 0.75 $ & 4.23\e{-15}   \\ & 
$0.75 - 1.00 $ & 4.37\e{-15} & & 
$0.75 - 1.00 $ & 4.42\e{-15}   \\\addlinespace[.25em]
$4^6$ to $4^7$                                                                                       &$0.00 - 0.25 $ & 5.01\e{-15} &   \\ & 
$0.25 - 0.50 $ & 3.63\e{-15} &   \\ & 
$0.50 - 0.75 $ & 3.97\e{-15} &   \\ & 
$0.75 - 1.00 $ & 4.48\e{-15} &   \\\addlinespace[.25em]
\bottomrule
\end{tabular}

\end{center}
\caption{The maximum relative error encountered while evaluating $\chi_n(\gamma)$ for various
ranges of the parameters.}
\label{experiments:table1}
\end{table}

\end{subsection}


\begin{subsection}{The time with required to calculate $\chi_\sigma(\gamma)$}

In the experiments described here, we measured the time required to evaluate $\chi_\sigma(\gamma)$
using the algorithm of this paper.  

Table~\ref{experiments:table2} reports the results of the first set of such experiments.
In each experiment,  we fixed  a range $[\gamma_1,\gamma_2]$ of values of $\gamma$ and a range $[\sigma_1,\sigma_2]$ 
of values of $\sigma$.  
We sampled $100$ random values of $\gamma$ in  $[\gamma_1,\gamma_2]$
and $100$ random integers $n$ in  the interval $[\sigma_1\gamma,\sigma_2\gamma]$.
We then measured the time required to evaluate $\chi_n(\gamma)$
at each of the $10,000$ pairs of the sampled parameters using our algorithm
and using the Osipov-Xiao-Rokhlin method.  The average time required
by each approach is reported in Table~\ref{experiments:table2}.

\begin{table}[!h]
\begin{center}
\small
\begin{tabular}{ccccc}
\toprule
Range of $\gamma$ & Range of $\sigma$ & Average Time   & Average time    \\                               &                   &  expansion     & Rokhlin, et. al.\\
\midrule
$4^3$ to $4^4$                                                                                       &$0.00 - 0.25 $ & $1.04\e{-06}$ &$5.28\e{-05}$ \\  \addlinespace[.125em]
 & 
$0.25 - 0.50 $ & $1.09\e{-06}$ &$6.86\e{-05}$ \\  \addlinespace[.125em]
 & 
$0.50 - 0.75 $ & $1.25\e{-06}$ &$8.40\e{-05}$ \\  \addlinespace[.125em]
 & 
$0.75 - 1.00 $ & $1.41\e{-06}$ &$1.01\e{-04}$ \\  \addlinespace[.125em]
 \addlinespace[.25em]
$4^4$ to $4^5$                                                                                       &$0.00 - 0.25 $ & $1.08\e{-06}$ &$1.38\e{-04}$ \\  \addlinespace[.125em]
 & 
$0.25 - 0.50 $ & $1.17\e{-06}$ &$2.16\e{-04}$ \\  \addlinespace[.125em]
 & 
$0.50 - 0.75 $ & $1.38\e{-06}$ &$2.97\e{-04}$ \\  \addlinespace[.125em]
 & 
$0.75 - 1.00 $ & $1.59\e{-06}$ &$3.64\e{-04}$ \\  \addlinespace[.125em]
 \addlinespace[.25em]
$4^5$ to $4^6$                                                                                       &$0.00 - 0.25 $ & $1.18\e{-06}$ &$4.17\e{-04}$ \\  \addlinespace[.125em]
 & 
$0.25 - 0.50 $ & $1.29\e{-06}$ &$8.30\e{-04}$ \\  \addlinespace[.125em]
 & 
$0.50 - 0.75 $ & $1.49\e{-06}$ &$1.06\e{-03}$ \\  \addlinespace[.125em]
 & 
$0.75 - 1.00 $ & $1.69\e{-06}$ &$1.41\e{-03}$ \\  \addlinespace[.125em]
 \addlinespace[.25em]
$4^6$ to $4^7$                                                                                       &$0.00 - 0.25 $ & $1.27\e{-06}$ &$1.62\e{-03}$ \\  \addlinespace[.125em]
 & 
$0.25 - 0.50 $ & $1.34\e{-06}$ &$3.01\e{-03}$ \\  \addlinespace[.125em]
 & 
$0.50 - 0.75 $ & $1.54\e{-06}$ &$4.26\e{-03}$ \\  \addlinespace[.125em]
 & 
$0.75 - 1.00 $ & $1.76\e{-06}$ &$5.29\e{-03}$ \\  \addlinespace[.125em]
 \addlinespace[.25em]
$4^7$ to $4^8$                                                                                       &$0.00 - 0.25 $ & $1.33\e{-06}$ &$6.40\e{-03}$ \\  \addlinespace[.125em]
 & 
$0.25 - 0.50 $ & $1.36\e{-06}$ &$1.21\e{-02}$ \\  \addlinespace[.125em]
 & 
$0.50 - 0.75 $ & $1.61\e{-06}$ &$1.73\e{-02}$ \\  \addlinespace[.125em]
 & 
$0.75 - 1.00 $ & $1.86\e{-06}$ &$2.13\e{-02}$ \\  \addlinespace[.125em]
 \addlinespace[.25em]
$4^8$ to $4^9$                                                                                       &$0.00 - 0.25 $ & $1.46\e{-06}$ &$2.50\e{-02}$ \\  \addlinespace[.125em]
 & 
$0.25 - 0.50 $ & $1.49\e{-06}$ &$4.60\e{-02}$ \\  \addlinespace[.125em]
 & 
$0.50 - 0.75 $ & $1.77\e{-06}$ &$6.73\e{-02}$ \\  \addlinespace[.125em]
 & 
$0.75 - 1.00 $ & $1.98\e{-06}$ &$7.56\e{-02}$ \\  \addlinespace[.125em]
 \addlinespace[.25em]
$4^9$ to $4^{10}$                                                                                    &$0.00 - 0.25 $ & $1.58\e{-06}$ &$9.86\e{-02}$ \\  \addlinespace[.125em]
 & 
$0.25 - 0.50 $ & $1.59\e{-06}$ &$1.86\e{-01}$ \\  \addlinespace[.125em]
 & 
$0.50 - 0.75 $ & $1.82\e{-06}$ &$2.81\e{-01}$ \\  \addlinespace[.125em]
 & 
$0.75 - 1.00 $ & $2.04\e{-06}$ &$3.22\e{-01}$ \\  \addlinespace[.125em]
 \addlinespace[.25em]
\bottomrule
\end{tabular}

\end{center}
\caption{The average time (in seconds) required to evaluate $\chi_n(\gamma)$ using the algorithm of this
paper and the Osipov-Xiao-Rokhlin algorithm for various ranges of the parameters.}
\label{experiments:table2}
\end{table}

In a second set of experiments, we measured the time required to evaluate $\chi_\sigma(\gamma)$
as $\gamma$ varies for certain fixed values of $\sigma$ and 
the time required to evaluate $\chi_\sigma(\gamma)$ as $\sigma$ varies for certain
fixed values of $\gamma$.    Figure~\ref{experiments:figure1}  gives the results.

\begin{figure}[!h]
\hfil
\includegraphics[width=.45\textwidth]{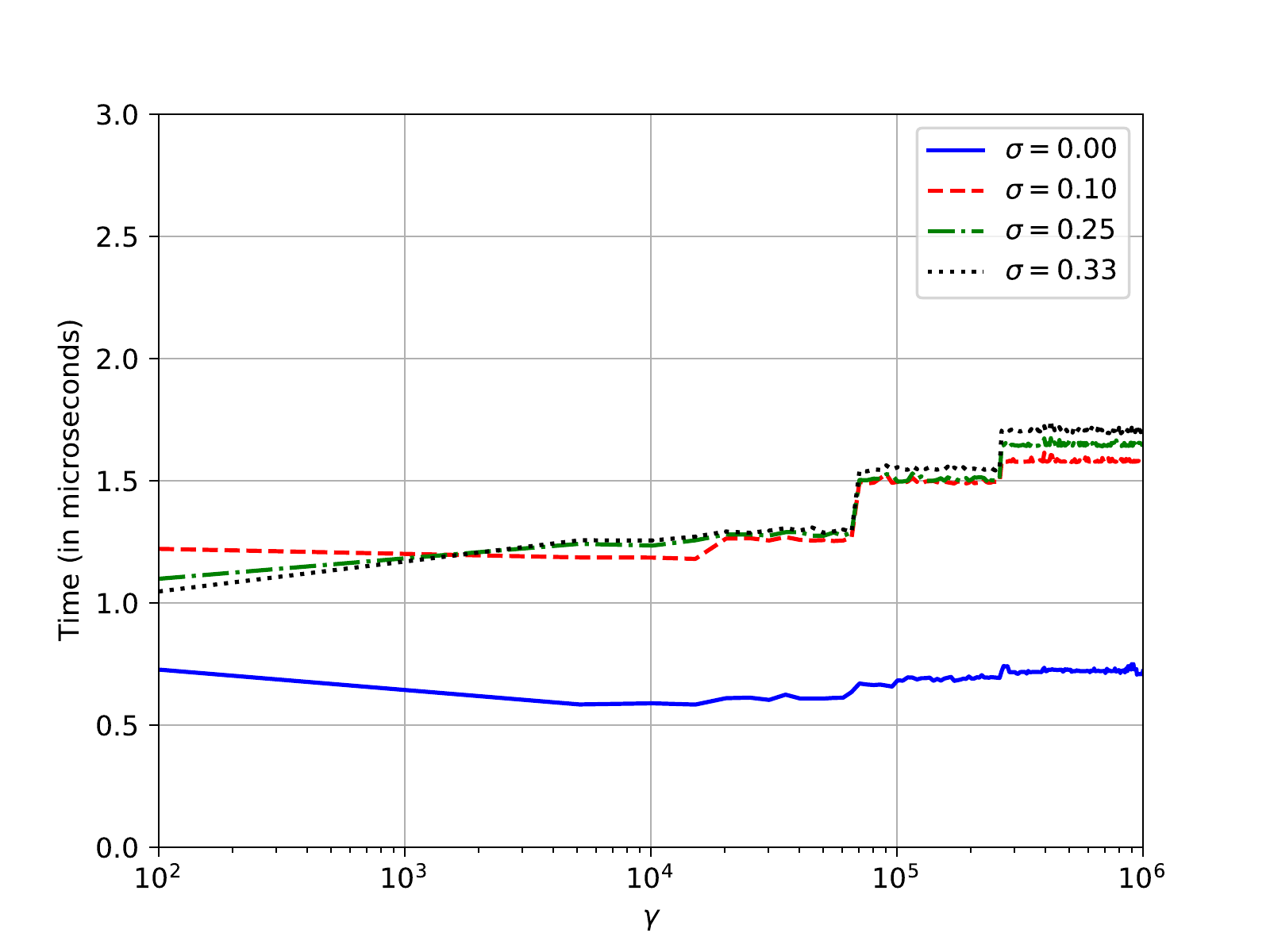}
\hfil
\includegraphics[width=.45\textwidth]{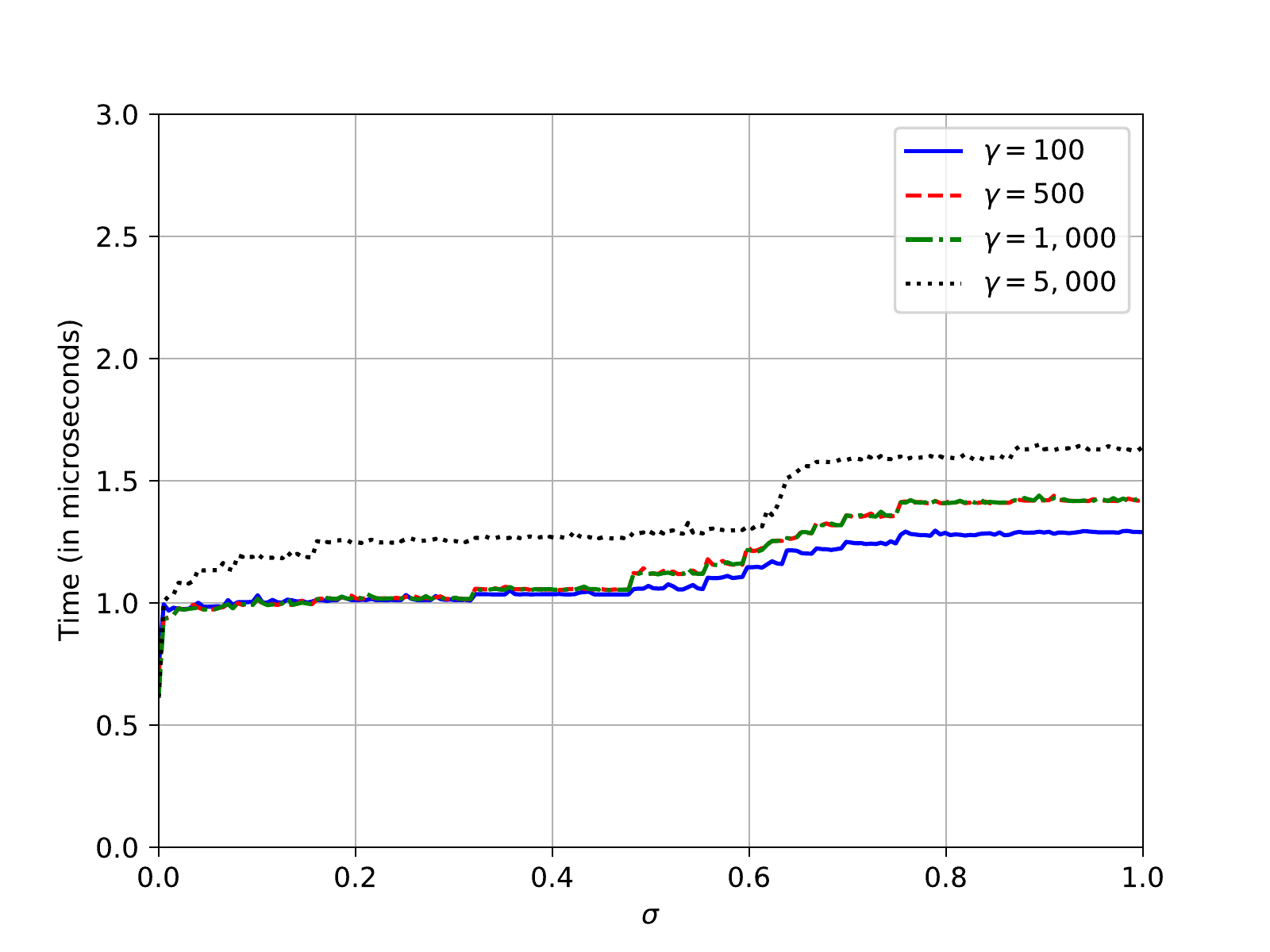}
\hfil

\hfil
\includegraphics[width=.45\textwidth]{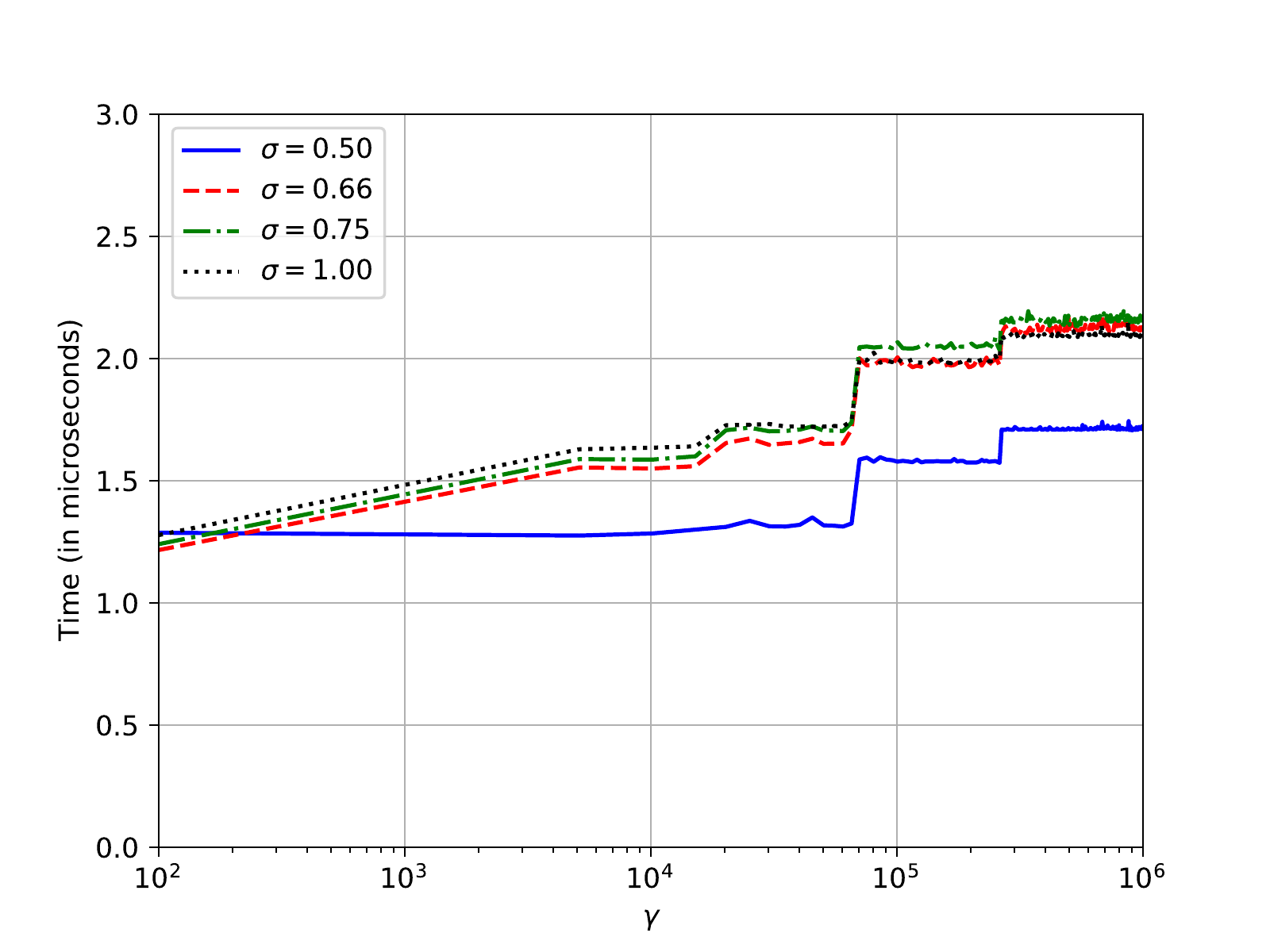}
\hfil
\includegraphics[width=.45\textwidth]{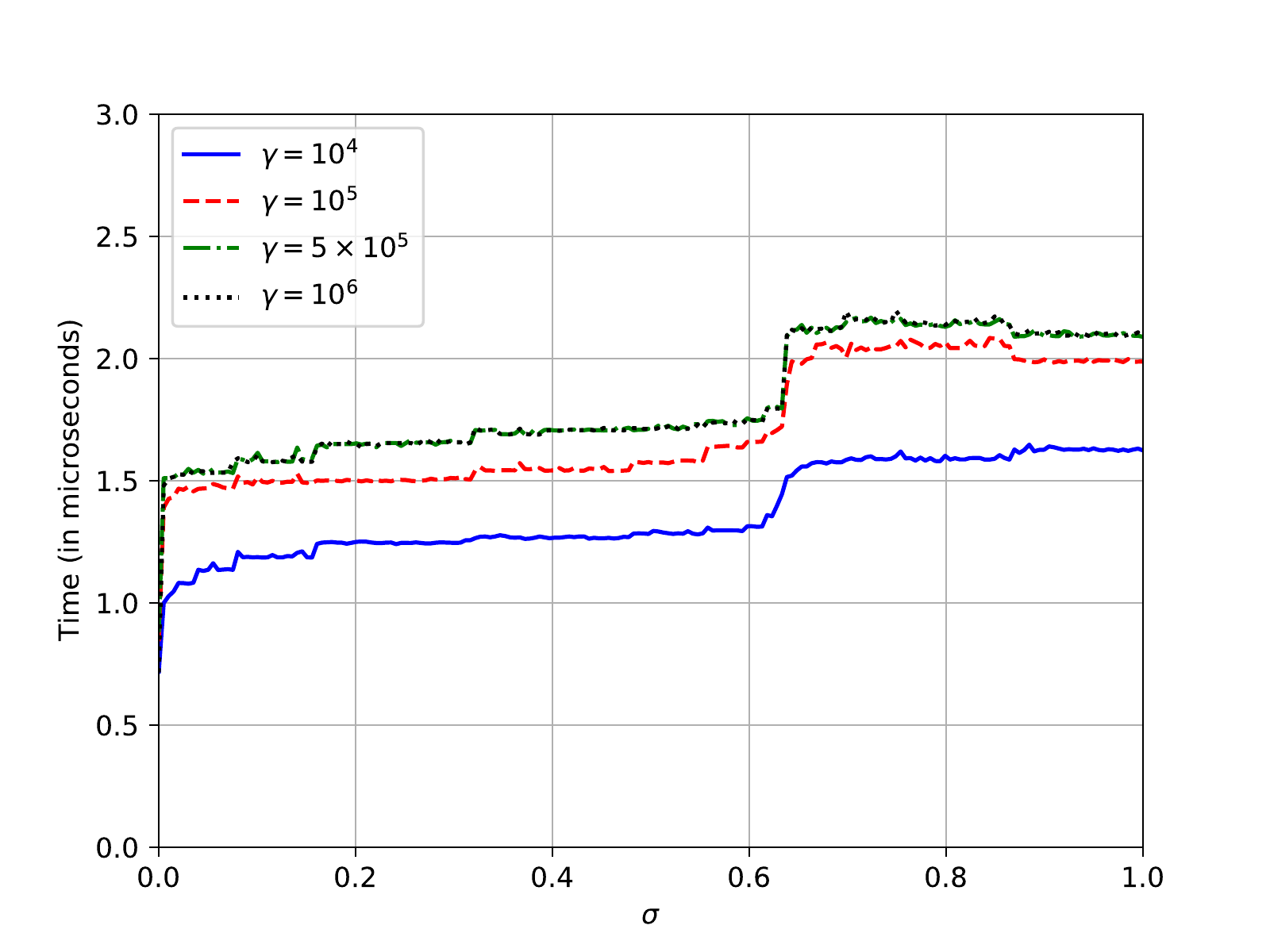}
\hfil

\caption{The time (in microseconds) required by the algorithm of this article to evaluate $\chi_\sigma(\gamma)$.
Each plot on the left gives the time as a function of $\gamma$ for various values of $\sigma$, 
and each plot on the right gives the time as a function of $\sigma$ for various values of $\gamma$.
A logarithmic scale is used for the x-axis in each of the graphs on the left.}

\label{experiments:figure1}
\end{figure}

\end{subsection}

\end{section}

\begin{section}{Acknowledgements}
The second author was supported in part by an NSERC Discovery grant  RGPIN-2021-02613, and
by NSF grants DMS-1818820 and  DMS-2012487.
\end{section}

\vfill\eject
\bibliographystyle{acm}
\bibliography{prolates.bib}

\end{document}